\documentclass[secthm]{elsart_N}

\usepackage{amssymb,amsmath}
\usepackage{graphics}
\usepackage[arrow, matrix, curve]{xy}

\usepackage{bbm}

\theoremstyle{plain}
\newtheorem{thmT}{Theorem}

\newenvironment{thmcite}[1]{
\textbf{Theorem (#1)}\it}

\newenvironment{nil}{\it }

\renewcommand{\Cup}{\bigcup}

\newcommand{\restr}{\upharpoonright}
\newcommand{\standPol}{Let $X,Y$ be Polish spaces equipped with Borel probability measures $\mu,\nu$}

\newcommand{\N}{\mathbb{N}}
\newcommand{\Q}{\mathbb{Q}}
\newcommand{\R}{\mathbb{R}}
\newcommand{\Rp}{(-\infty, \infty]}
\newcommand{\Rm}{[-\infty, \infty)}
\newcommand{\Rop}{[0,\infty]}

\newcommand{\eps}{\varepsilon}
\renewcommand{\phi}{\varphi}
\renewcommand{\rho}{\varrho}
\newcommand{\me}{^{-1}}
\newcommand{\push}{{}_{\#}}
\newcommand{\Eins}{\mathbbm{1}}
\renewcommand{\labelenumi}{\alph{enumi}.}

\begin{document}
\begin{frontmatter}

\title{Optimal and Better Transport Plans}
\author[UNI]{Mathias Beiglb\"ock\thanksref{S9612}},
\author[DMG]{Martin Goldstern},
\author[DMG]{Gabriel Maresch \thanksref{Y328}},
\author[UNI]{Walter Schachermayer\thanksref{multi}}
\thanks[S9612]{Supported by the Austrian Science Fund (FWF) under grant S9612.}
\thanks[Y328]{Supported by the Austrian Science Fund (FWF) under grant Y328 and P18308.}
\thanks[multi]{Supported by the Austrian Science Fund (FWF) under grant
P19456, from the Vienna Science and Technology Fund (WWTF) under
grant MA13 and from the Christian Doppler Research Association
(CDG)}

\address[DMG]{Institut f\"ur Diskrete Mathematik und Geometrie, Technische Universit\"at Wien\endgraf
Wiedner Hauptstra\ss e 8-10/104,1040 Wien, Austria}

\address[UNI]{Fakult\"at f\"ur Mathematik, Universit\"at Wien\endgraf
Nordbergstra\ss e 15, 1090 Wien, Austria}


\begin{abstract}
We consider the Monge-Kantorovich transport problem in a purely measure theoretic setting, i.e.\ without imposing continuity assumptions on the cost function. It is known that transport plans which are concentrated on $c$-monotone sets are optimal, provided the cost function $c$ is either lower
semi-continuous and finite, or continuous and may possibly attain the value $\infty$. We
show that this is true in a more general setting, in particular for merely Borel measurable  cost functions provided that $\{c=\infty\}$ is the union of a closed set and a negligible set.
In a previous paper Schachermayer and Teichmann considered strongly $c$-monotone transport plans and proved that every strongly $c$-monotone transport plan is optimal. We establish that transport plans are strongly $c$-monotone if and only if they satisfy a ``better'' notion of optimality called robust optimality.
\end{abstract}
\end{frontmatter}

\section{Introduction}\label{Sintro}

We consider the \emph{Monge-Kantorovich transport problem}
$(\mu,\nu,c)$ for Borel probability  measures $\mu,\nu$ on Polish
spaces $X,Y$ and a Borel measurable cost function $c:X\times Y\to
\Rop$. As standard references on the theory of mass transport we
mention \cite{AmPr03,RaRu98,Vill03,Vill05}. By $\Pi(\mu,\nu)$ we
denote the set of all probability measures on $X\times Y$ with
$X$-marginal $\mu$ and $Y$-marginal $\nu$. For a Borel  measurable
\emph{cost function} $c:X\times Y\to\Rop$ the transport costs of a
given transport plan $\pi\in\Pi(\mu,\nu)$ are defined by
\begin{equation}\label{eqncost}
I_c[\pi]:=\int_{X\times Y} c(x,y) d\pi.
\end{equation} $\pi$ is called a \emph{finite} transport plan if $I_c[\pi]<\infty$.

A nice  interpretation of the Monge-Kantorovich transport problem is given by C{\'e}dric Villani in Chapter 3 of the impressive monograph \cite{Vill05}:

\medskip\begin{quotation}\begin{small}
``Consider a large number of bakeries, producing breads, that should be transported each morning to caf{\'e}s where consumers will eat them. The amount of bread that can be produced at each bakery, and the amount that will be consumed at each caf{\'e} are known in advance, and can be modeled as probability measures (there is a ``density of production'' and a ``density of consumption'') on a certain space, which in our case would be Paris (equipped with the natural metric such that the distance between two points is the length of the shortest path joining them). The problem is to find in practice where each unit of bread should go, in such a way as to minimize the total transport cost."\end{small}\end{quotation}\medskip

We are interested in \emph{optimal} transport plans, i.e.\
minimizers of the functional $I_c[\cdot]$
and their characterization via the notion of $c$-monotonicity.
\begin{defn}\label{Dcmon}
A Borel set $\Gamma\subseteq X\times Y$ is called
\emph{$c$-monotone} if
\begin{equation}\label{eqncmon}
\sum_{i=1}^n c(x_i,y_i)\leq \sum_{i=1}^n c(x_{i},y_{i+1})
\end{equation}
for all pairs $(x_1,y_1),\dotsc, (x_n,y_n)\in \Gamma$ using the
convention $y_{n+1}:=y_1$.  A transport plan $\pi$ is called
$c$-monotone if there exists a $c$-monotone $\Gamma$ with $\pi(\Gamma)=1$.
\end{defn}
In the literature (e.g.\ \cite{AmPr03,GaMc96,KnSm92,Prat07,ScTe08}) the following characterization was established under various continuity assumptions on the cost function.
Our main result states that those assumptions are not required.

\begin{thmT}\label{CmonEquivOptimal}\standPol\ and let $c: X\times Y \to \Rop$ a Borel measurable cost function.

\begin{enumerate}
\item \label{opt->cmon} Every finite optimal transport plan is $c$-monotone.
\item \label{cmon->opt} Every finite $c$-monotone
 transport plan is optimal
 if there exist a closed set $F$ and a $\mu\otimes\nu$-null set $N$
  such that $\{(x,y): c(x,y)=\infty\} = F\cup N$.
\end{enumerate}
\end{thmT}
Thus in the case of a cost function which does not attain the value $\infty$ the equivalence of optimality and $c$-monotonicity is valid without any restrictions beyond the obvious measurability conditions inherent in the formulation of the problem.

The subsequent construction due to Ambrosio and Pratelli in
\cite[Example 3.5]{AmPr03} shows that if $c$ is allowed to attain
$\infty$ the implication ``$c$-monotone $\Rightarrow$ optimal'' does not hold without some additional assumption as in Theorem \ref{CmonEquivOptimal}.b.


\begin{exmp}[Ambrosio and Pratelli]\label{ExAmPra}
Let $X=Y= [0,1]$, equipped with Lebesgue measure $\lambda=\mu=\nu$.
Pick $\alpha\in [0,1)$ irrational. Set $$\Gamma_0=\{(x,x):x\in
X\},\quad\Gamma_1=\{(x,x\oplus\alpha):x\in X\},$$ where $\oplus$ is
addition modulo $1$. Let  $c: X\times Y \to \Rop$ be such that
$c=a\in[0,\infty)$ on $\Gamma_0$, $c=b\in[0,\infty)$ on $ \Gamma_1$
and $c=\infty$ otherwise. It is then easy to check that  $\Gamma_0$
and $\Gamma_1$ are $c$-monotone sets. Using the maps $f_0, f_1: X
\to X\times Y$, $f_0(x)=(x,x), f_1(x) = (x,x\oplus \alpha)$ one
defines the transport plans $\pi_0=f_0\push\lambda$, $\pi_1=f_1\push
\lambda$ supported by $\Gamma_0$ respectively $\Gamma_1$. Then
$\pi_0$ and $\pi_1$ are finite $c$-monotone transport plans, but as
$I_c[\pi_0]=a, I_c[\pi_1]=b$ it depends on the choice of $a$ and $b$
which transport plan is optimal. Note that in contrast to the
assumption in Theorem \ref{CmonEquivOptimal}.b the set $\{(x,y)\in
X\times Y:c=\infty\}$ is open.
\end{exmp}

We want to remark that rather trivial (folkloristic) examples show
that no optimal transport has to exist if the cost function doesn't
satisfy proper continuity assumptions.
\begin{exmp} Consider the task to transport points on
the real line (equipped  with the Lebesgue measure) from the
interval $[0,1)$ to $[1,2)$ where the cost of moving one point to
another is the squared distance between these points ($X=[0,1), Y=
[1,2)$, $c(x,y)=(x-y)^2$, $\mu=\nu=\lambda$). The simplest way to
achieve this transport is to shift every point by $1$. This results
in transport costs of $1$ and one easily checks that all other
transport plans are more expensive.

If we now alter the cost function to be $2$ whenever two points have
distance $1$, i.e.\ if we set $$\tilde c(x,y)=
\begin{cases} 2& \mbox{if } y=x+1 \\ c(x,y) & \mbox{otherwise} \end{cases},$$
it becomes impossible to find a transport plan $\pi\in\Pi(\mu,\nu)$
with total transport costs $I_{\tilde c} [\pi]=1$, but it is still
possible to achieve transport costs arbitrarily close to $1$. (For
instance, shift $[0,1-\eps)$ to $[1+\eps,2)$ and $[1-\eps,1)$ to
$[1,1+\eps)$ for small $\eps>0$.)
\end{exmp}

\subsection{History of the problem}
The notion of $c$-monotonicity originates in convex analysis. The
well known Rockafellar Theorem (see for instance \cite[Theorem
3]{Rock66} or \cite[Theorem 2.27]{Vill03}) and its generalization,
R\"uschendorf's Theorem (see \cite[Lemma 2.1]{R96}), characterize
$c$-monotonicity in $\R^n$ in terms of integrability. The
definitions of $c$-concave functions and super-differentials can be
found for instance in \cite[Section 2.4]{Vill03}.

\begin{thmcite}{Rockafellar}\label{Trockafellar}
A non-empty set $\Gamma\subseteq \R^n\times \R^n$ is  cyclically monotone (that is,  $c$-monotone with respect to the squared euclidean distance) if and only if there exists a l.s.c.\ concave function $\phi: \R^n \to \R$ such that $\Gamma$ is contained in the super-differential $\partial(\phi)$.
\end{thmcite}

\begin{thmcite}{R\"uschendorf}\label{Trueschendorf}
Let $X,Y$ be abstract spaces and $c: X\times Y\to \Rop$ arbitrary. Let   $\Gamma\subseteq X\times Y$ be $c$-monotone. Then there exists a $c$-concave function $\phi: X \to Y$ such that $\Gamma$ is contained in the $c$-super-differential $\partial^c(\phi)$.
\end{thmcite}

Important results of Gangbo and McCann \cite{GaMc96} and Brenier
\cite[Theorem 2.12]{Vill03} use these potentials to establish
uniqueness of the solutions of the Monge-Kantorovich transport
problem in $\R^n$ for different types of cost functions subject to
certain regularity conditions.

\medskip
\emph{Optimality implies $c$-monotonicity: } This is evident in the
discrete case if $X$ and $Y$ are finite sets. For suppose that $\pi$
is a transport plan for which $c$-monotonicity is violated on pairs
$(x_1,y_1),\dotsc,(x_n,y_n)$ where all points $x_1,\dotsc, x_n$ and
$y_1,\dotsc, y_n$ carry positive mass. Then we can reduce costs by
sending the mass $\alpha>0$, for $\alpha$ sufficiently small, from
$x_i$ to $y_{i+1}$ instead of $y_i$, that is, we replace the
original transport plan $\pi$ with
\begin{equation}\pi^\beta=\pi+\alpha
\sum_{i=1}^n\delta_{(x_i,y_{i+1})} -\alpha
\sum_{i=1}^n\delta_{(x_i,y_i)}.\end{equation} (Here we are using the
convention $y_{n+1}=y_1.$)

 Gangbo and McCann (\cite[Theorem 2.3]{GaMc96}) show how continuity assumptions on the cost function can be exploited to extend this to an abstract setting. Hence one achieves:

\begin{nil}\label{Popt->cmon}
Let $X$ and $Y$ be Polish spaces equipped with Borel probability
measures $\mu,\nu$. Let $c: X\times Y\to \Rop$ be a l.s.c.\ cost
function. Then every finite optimal transport plan is $c$-monotone.
\end{nil}

Using measure theoretic tools, as developed in the beautiful paper by Kellerer \cite{Kell84}, we are able to extend this to Borel measurable cost functions (Theorem \ref{CmonEquivOptimal}.a.) without any additional regularity assumption. 

\medskip

\emph{$c$-monotonicity implies optimality:} In the case of finite
spaces $X,Y$ this again is nothing more than an easy exercise
(\cite[Exercise 2.21]{Vill03}). The problem gets harder in the
infinite setting. It was first proved in \cite{GaMc96} that for
$X,Y$ compact subsets of $\R^n$ and $c$ a continuous cost function,
$c$-monotonicity implies optimality. In a more general setting this
was shown in \cite[Theorem 3.2]{AmPr03} for l.s.c.\ cost functions
which additionally satisfy the moment conditions
\begin{align*}
\mu\;\Bigg(\Bigg\{ x: \int_Y c(x,y) d\nu < \infty\Bigg\}\Bigg) &> 0,\\
\nu\;\Bigg(\Bigg\{ y: \int_X c(x,y) d\mu < \infty\Bigg\}\Bigg) &> 0.
\end{align*}

Further research into this direction was initiated by the following problem posed by Villani in \cite[Problem 2.25]{Vill03}:

\begin{nil} For $X=Y=\R^n$ and $c(x,y)=\|x-y\|^2$, the squared euclidean distance, does $c$-monotonicity of a transport plan imply its optimality?\end{nil}
A positive answer to this question was given independently by
Pratelli in \cite{Prat07} and by Schachermayer and Teichmann in
\cite{ScTe08}. Pratelli proves the result for countable spaces and
shows that it  extends to the  Polish case by means of
approximation if the cost function $c:X\times Y\to\Rop$ is
continuous. The paper \cite{ScTe08} pursues a different approach:
The notion of \emph{strong $c$-monotonicity} is introduced. From
this property optimality follows fairly easily and the main part of
the paper is concerned with the fact that strong $c$-monotonicity
follows from the usual notion of $c$-monotonicity in the Polish
setting if $c$ is assumed to be l.s.c.\ and finitely valued.

Part (b) of Theorem \ref{CmonEquivOptimal} unifies these statements:
Pratelli's result follows from the fact that for continuous $c:
X\times Y \to \Rop$ the set $\{ c = \infty \}=c^{-1}[\{\infty\}]$ is
closed; the Schachermayer-Teichmann result follows since for finite
$c$ the set $\{ c = \infty \}$ is empty.

Similar to \cite{ScTe08} our proofs are based on the concept of
\emph{strong $c$-monotonicity}. In Section \ref{StrongNotions} we
present \emph{robust optimality} which is a variant of optimality
that we shall show to be equivalent to strong $c$-monotonicity. As
not every optimal transport plan is also robustly optimal, this
accounts for the somewhat provocative concept of ``better than
optimal'' transport plans alluded to in the title of this paper.

Correspondingly the notion  of strong $c$-monotonicity is in fact
stronger than ordinary $c$-monotonicity (at least if $c$ is allowed
to assume the value $\infty$). 

\subsection{Strong Notions}\label{StrongNotions}

It turns out that optimality of a transport plan  is intimately
connected with the notion of  \emph{strong $c$-monotonicity}
introduced in \cite{ScTe08}.
\begin{defn}\label{Dstrong} A Borel set $\Gamma\subseteq X\times Y$
 is \emph{strongly $c$-monotone} if
   there exist Borel measurable functions
  $\phi:X\to \Rm$ and $\psi:Y\to\Rm$ such that $\phi(x)+\psi(y)\leq c(x,y)$ for all $(x,y)\in X\times Y$ and $ \phi(x)+ \psi(y)=c(x,y)$ for all $(x,y)\in\Gamma$.
A transport plan $\pi\in \Pi(\mu,\nu)$ is \emph{strongly
$c$-monotone} if $\pi$ is concentrated on  a strongly $c$-monotone
Borel set $\Gamma$.
\end{defn}
Strong $c$-monotonicity implies $c$-monotonicity since
\begin{align}\label{eqnstrong}
\sum_{i=1}^n c(x_{i+1},y_{i})\ge\sum_{i=1}^n \phi(x_{i+1})+\psi(y_{i})=
\sum_{i=1}^n \phi(x_i)+\psi(y_i)=\sum_{i=1}^n c(x_i,y_i)
\end{align}
whenever $(x_1,y_1),\dotsc, (x_n,y_n)\in \Gamma$ .

If there are \emph{integrable} functions $\phi$ and $\psi$ witnessing that $\pi$ is strongly $c$-monotone, then for every $\tilde \pi\in\Pi(\mu,\nu)$ we can estimate:
\begin{align*}
I_c[\pi] =& \int_{\Gamma} c(x,y)d\pi=\int_{\Gamma} [\phi(x)+\psi(y)]d\pi=\\
         =&\int_{\Gamma} \phi(x)d\mu +\int_{\Gamma}
\psi(y)d\nu=\int_{\Gamma} [\phi(x)+\psi(y)]d\tilde\pi \le
I_c[\tilde\pi].
\end{align*}
Thus in this case strong $c$-monotonicity implies optimality.
However there is no reason why the Borel measurable functions
$\phi,\psi$ appearing in Definition \ref{Dstrong} should be integrable. In \cite[Proposition 2.1]{ScTe08} it is shown that for
l.s.c.\ cost functions, there is a way of truncating which allows to
also handle non-integrable functions $\phi$ and $\psi$.
The proof extends to merely Borel measurable functions; hence we have:
\begin{prop}\label{StrongCmonImpliesOptimal}
\standPol\ and
let $c: X\times Y\to \Rop$ be Borel measurable. Then every finite
transport plan which is strongly $c$-monotone is optimal.
\end{prop}
No new ideas are required to extend \cite[Proposition 2.1]{ScTe08} to the present setting but since Proposition \ref{StrongCmonImpliesOptimal} is a crucial ingredient of several proofs in this paper we provide an outline of the argument in Section 3.

As it will turn out, strongly $c$-monotone transport plans even
satisfy a ``better'' notion of optimality, called \emph{robust
optimality}.
\begin{defn}\label{Drobust} \standPol\ and let
$c:X\times Y\to \Rop$ be a Borel measurable cost function. A
transport plan $\pi\in \Pi(\mu,\nu)$ is \emph{robustly optimal} if,
for any Polish space $Z$ and any finite Borel measure $\lambda\ge 0$ on $Z$, there exists a Borel measurable extension $\tilde c: (X\cup
Z)\times (Y\cup Z)\to \Rop$ satisfying
\[\tilde{c}(a,b)=\left\{\begin{array}{ccl}
c(a,b)&\mbox{ for }& a\in X, b\in Y\\
0     &\mbox{ for }& a,b \in Z\\
<\infty&\mbox{ otherwise }\end{array} \right.\] such that the
measure
$\tilde\pi:=\pi+\left(\mbox{id}_Z\!\times\!\mbox{id}_Z\right)\push\lambda$
is optimal on $(X\cup Z)\times (Y\cup Z)$.  Note that $\tilde \pi$
is not a probability measure, but has total mass $1+\lambda(Z)\in
[1,\infty)$.
\end{defn}
Note that since we allow the possibility $\lambda (Z)=0 $ every robustly optimal transport plan is in particular optimal in the usual sense.

Robust optimality has a colorful ``economic'' interpretation: a tycoon wants to enter the Parisian croissant consortium. She builds a
storage of size $\lambda(Z)$ where she buys up croissants and sends
them to the caf\'{e}s. Her hope is that by offering low transport costs, the previously optimal
transport plan $\pi$ will not be optimal anymore, so that the
traditional relations between bakeries and caf\'{e}s will collapse.
Of course, the authorities of Paris will try to defend their
structure by imposing (possibly very high, but still finite) tolls
for all transports to and from the tycoon's storage, thus resulting in finite costs $\tilde c(a,b)$ for $(a,b)\in (X\times Z)\cup (Z\times Y)$. In the case of robustly optimal $\pi$ they can successfully defend themselves against the intruder.

Every robustly optimal transport $\pi$ plan is optimal in
the usual sense and hence also $c$-monotone. The crucial feature is
that robust optimality implies strong $c$-monotonicity. In fact, the
two properties are equivalent.
\begin{thmT}\label{StrongCmonEquivRobustlyOptimal}
\standPol\ and $c:X\times Y\to \Rop$ a
Borel measurable cost function. For a finite transport plan $\pi$
the following assertions are equivalent:
\begin{enumerate}
\item\label{T2strong} $\pi$ is strongly $c$-monotone.
\item\label{T2robust} $\pi$ is robustly optimal.
\end{enumerate}
\end{thmT}
Example \ref{ZeroOneExample} below shows that robust optimality resp.\
strong c-monotonicity is in fact a stronger property than usual
optimality. 

\subsection{Putting things together}

Finally we want to point out that  in the situation where $c$ is finite all previously mentioned  notions of monotonicity and optimality coincide. We can even pass to a slightly more general setting than finite cost functions and obtain the following result.
\begin{thmT}\label{AllEquiv}
\standPol\ and let $c: X\times Y\to \Rop$ be Borel measurable and
$\mu\otimes\nu$-a.e.\ finite. For a finite transport plan $\pi$ the
following assertions are equivalent:

\begin{enumerate}\renewcommand{\labelenumi}{(\arabic{enumi})}
\item \label{T3optimal} $\pi$ is optimal.
\item \label{T3cmon} $\pi$ is $c$-monotone.
\item \label{T3robust} $\pi$ is robustly optimal.
\item \label{T3strong} $\pi$ is strongly $c$-monotone.
\end{enumerate}
\end{thmT}

The equivalence of (1), (2) and (4) was established in \cite{ScTe08} under the additional assumption that $c$ is l.s.c.\ and finitely valued.
\bigskip

We sum up the situation under fully general assumptions.
The upper line (1 and 2)  relates to the optimality of a transport plan $\pi$. The lower line (3 and 4) contains the two equivalent strong concepts and implies the upper line but - without additional assumptions - not vice versa.

\begin{center}\begin{figure}[h]\label{diagram}
\fbox{
\begin{minipage}[r]{10cm}\parbox[c]{5cm}{
\begin{xy}
  \xymatrix{
   1 \ar@<-2pt>[r]   
   &2\ar@{.>}@<-2pt>[l]_{\text{Thm. 1}}
   \ar@<2pt>@{.>}[d]^{\text{\;\,Thm. 3}}
   &\\
    3\ar[u]\ar@{<->}[r]\ar@{}\ar@{}@<-2pt>[r]_{\text{Thm. 2}} &4&
}
\end{xy}}
\hfill
\parbox[b]{6cm}{
\begin{tabular}{cl}
(1)&optimal\\
(2)&$c$-monotone\\
(3)&robustly optimal\\
(4)&strongly $c$-monotone
\end{tabular}}
\end{minipage}}
\caption{Implications between properties of transport plans}
\end{figure}
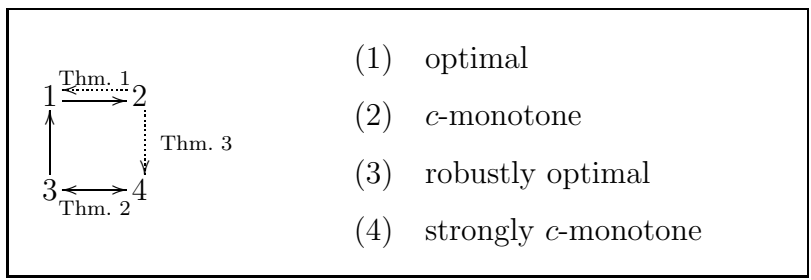
\end{center}
Note that the implications symbolized by dotted lines in Figure 1
are not true without additional assumptions ($(2) \not\Rightarrow
(1)$: Example \ref{ExAmPra}, $(1) \not\Rightarrow (3)$ resp.\ $(4)$: Example \ref{ZeroOneExample}).

The paper is organized as follows:
In Section 2 we prove that every optimal transport plan $\pi$ is $c$-monotone (Theorem \ref{CmonEquivOptimal}.a). In Section 3 we introduce an auxilliary property (\emph{connectedness}) of the support of a transport plan  and show that it allows to pass from $c$-monotonicity to strong $c$-monotonicity. Moreover we establish that strong $c$-monotonicity implies optimality (Proposition \ref{StrongCmonImpliesOptimal}). Section 4 is concerned with the proof of Theorem \ref{CmonEquivOptimal}.b. 
Finally we complete the proofs of Theorems \ref{StrongCmonEquivRobustlyOptimal} and \ref{AllEquiv} in Section 5.




We observe that in all the above discussion we only referred to the Borel structure of the Polish spaces $X,Y$, and never referred to the topological structure. Hence the above results (with the exception of Theorem \ref{CmonEquivOptimal}.b.) hold true for standard Borel measure spaces.

In fact it seems likely that our results can be transferred to the
setting of perfect measure spaces. (See \cite{Rama02} for a general
overview resp.\ \cite{RaRu98} for a treatment of problems of mass
transport in this framework.) However we do not pursue this
direction.

\noindent {\bf Acknowledgement.} The authors are indebted to the extremely careful referee who noticed many inaccuracies  resp.\ mistakes and whose insightful suggestions led to a more accessible presentation of several results in this paper.

\section{Improving Transports}\label{SmilestoneA}

Assume that some transport plan $\pi\in\Pi(\mu,\nu)$ is given. From
a purely heuristic point of view there are either
few tupels $\left((x_1,y_1),\dotsc, (x_n,y_n)\right)$
along which $c$-monotonicity is violated, or there are many such tuples, in which
case $\pi$ can be enhanced by rerouting the transport along these tuples.
As the notion of $c$-monotonicity
refers to $n$-tuples it turns out that it is necessary to
consider finitely many measure spaces to properly formulate what is
meant by ``few'' resp.\ ``many''.

 Let $X_1,\dotsc,X_n$ be
Polish spaces equipped with finite Borel measures
$\mu_1,\dotsc,\mu_n$. By $\Pi(\mu_1,\dotsc,\mu_n)\subseteq
\mathcal{M}(X_1\times\dotsb\times X_n)$ we denote the set of all
Borel measures on $X_1\times\dotsb\times X_n$ such that the $i$-th
marginal measure coincides with the Borel measure $\mu_i$ for
$i=1,\dotsc, n$. By $p_{X_i}: X_1\times\dotsb\times X_n\to X_i$ we
denote the projection onto the $i$-th component. $B\subseteq
X_1\times\dotsb\times X_n$ is called an \emph{L-shaped null set} if
there exist null sets $N_1\subseteq X_1,\dotsc,N_n\subseteq X_n$
such that $B\subseteq \bigcup_{i=1}^n p_{X_i}^{-1} [N_i]$.

The Borel sets of $X_1\times\dotsb\times X_n$ satisfy a nice
dichotomy. They are either L-shaped null sets or they carry a
positive measure whose marginals are absolutely continuous with
respect to $\mu_1,\dotsc,\mu_n$:
\begin{prop}\label{PKmeasure}
 Let $X_1,\ldots,X_n,n\geq 2$ be Polish spaces equipped with Borel probability measures $\mu_1,\dotsc,\mu_n$. Then for any Borel set  $B\subseteq X_1\times\dotsb\times X_n$ let
\begin{align}
P(B):=&\sup\left\{ \pi(B): \pi\in\Pi(\mu_1,\dotsc,\mu_n)\right\}\\
L(B):=&\inf\left\{\sum_{i=1}^n \mu_i(B_i): B_i \subseteq X_i \mbox{ and } B\subseteq \bigcup_{i=1}^n p_{X_i}^{-1}[B_i]\right\}.
\end{align}
Then $P(B)\geq1/n\,L(B)$. In particular $B$ satisfies one of the following alternatives:
\begin{enumerate}
\item $B$ is an L-shaped null set.
\item There exists $\pi\in \Pi(\mu_1,\dotsc,\mu_n)$ such that $\pi(B)>0$.
\end{enumerate}
\end{prop}

The main ingredient in the proof Proposition \ref{PKmeasure} is the following duality theorem due to Kellerer (see \cite[Lemma 1.8(a), Corollary 2.18]{Kell84}).

\begin{thmcite}{Kellerer}\label{KellererDuality}
  Let $X_1,\ldots,X_n,n\geq 2$ be Polish spaces equipped with Borel probability measures $\mu_1,\dotsc, \mu_n$ and assume that $c: X=X_1\times\dotsb\times X_n\to\R$ is Borel measurable and that
  $\overline c := \sup_X c, \underline c:=\inf_X c$ are finite.
  Set \begin{align*}
I(c)=&\inf\left\{ \int_X c\ d\pi: \pi\in\Pi(\mu_1,\dotsc,\mu_n)\right\},\\
S(c)=&\sup\left\{\sum_{i=1}^n\int_{X_i} \phi_i \ d\mu_i:
c(x_1,\dotsc,x_n)\geq \sum_{i=1}^n \phi_i(x_i), \frac {1}{n}
\overline c -(\overline c-\underline c)\leq \phi_i\leq \frac {1} {n}
\overline c \right\}.
\end{align*}
Then $I(c)=S(c)$.
\end{thmcite}
\begin{pf*}{PROOF of Proposition \ref{PKmeasure}.}
Observe that $-I(-\Eins_B)=P (B)$ and that
\begin{equation}\label{reversedDefinition} -S(-\Eins_B)=
\inf\!\left\{ \sum_{i=1}^n \int_{X_i}\chi_i d\mu_i:\Eins_B(x_1,\dotsc,
x_n)\le\sum_{i=1}^n \chi_i (x_i), 0\le \chi_i\le
1\right\}\!.\end{equation}
By Kellerer's Theorem 
 $-S(-\Eins_B)=-I(-\Eins_B)$. Thus it remains to show that $-S(-\Eins_B)\geq 1/n\,  L (B)$.
Fix functions $\chi_1,\ldots,\chi_n$ as in
(\ref{reversedDefinition}). Then  for each $(x_1,\ldots,x_n)\in B$
one has $1=\Eins_B(x_1,\ldots, x_n)\le\sum_{i=1}^n \chi_i (x_i)$ and
hence there exists some $i$ such that $\chi_i(x_i)\geq 1/n$. Thus
$B\subseteq\Cup_{i=1}^n p_{X_i}^{-1}[\{\chi_i\geq 1/n\}]$. It follows
that
\begin{eqnarray*}
-S(-\Eins_B)&\geq& \inf\left\{ \sum_{i=1}^n \int_{X_i}\chi_i d\mu_i:
 B\subseteq \Cup_{i=1}^n p_{X_i}^{-1}[\{\chi_i\geq 1/n\}], 0\le \chi_i\le 1\right\}\\
 &\geq& \inf\left\{ \sum_{i=1}^n \frac1n \mu_i(\{\chi_i\geq 1/ n\}) :
   B\subseteq \Cup_{i=1}^n p_{X_i}^{-1}[\{\chi_i\geq 1/n\}] \right\}\geq \frac1n L(B)
\end{eqnarray*} From this we deduce that either $L(B) =0$ or that there exists $\pi\in \Pi(\mu_1,\dotsc,\mu_n)$ such that $\pi(B)>0$.
The last assertion of Proposition \ref{PKmeasure} now follows from
the following Lemma due to Rich\'ard Balka and M\'arton Elekes
(private communication). \qed\end{pf*}
\begin{lem}\label{Lelekes}
Suppose that $L(B)=0$  for  a Borel set $B\subseteq X_1\times\dotsb\times X_n$. Then $B$ is an L-shaped null set.
\end{lem}

\begin{pf}
Fix $\eps>0$ and Borel sets $B_1^{(k)},\dotsc, B_n^{(k)}$ with
$\mu_i(B_i^{(k)})\le \eps\, 2^{-k}$ such that for each $k$
\[B\subseteq p_{X_1}^{-1}[B_1^{(k)}]\cup\dotsc\cup p_{X_n}^{-1}[B_n^{(k)}].\]
Let $B_i:=\bigcup_{k=1}^{\infty} B_i^{(k)}$ for $i=2,\dotsc,n$ such that
\[B\subseteq p_{X_1}^{-1}[B_1^{(k)}]\cup p_{X_2}^{-1}[B_2]\cup\dotsc\cup p_{X_n}^{-1}[B_n]\]
for each $k\in \N$. Thus with $B_1:=\bigcap_{k=1}^{\infty}
B_1^{(k)}$, \[B\subseteq p_{X_1}^{-1}[B_1]\cup
p_{X_2}^{-1}[B_2]\cup\dotsc\cup p_{X_n}^{-1}[B_n].\] Hence we can
assume from now on that $\mu_1(B_1)=0$ and that $\mu_i(B_i)$ is
arbitrarily small for $i=2,\dotsc, n$. Iterating this argument in
the obvious way we get the statement. \qed\end{pf}

\begin{rem}
In the case  $n=2$ it was shown in \cite[Proposition 3.3]{Kell84} that $L(B)=P(B)$ for every Borel set $B\subseteq X_1\times X_2$. However, for $n>2$, equality does not hold true, cf.\ \cite[Example 3.4]{Kell84}.
\end{rem}

\begin{defn}\label{Dbneps}
Let $X,Y$ be Polish spaces. For a Borel measurable cost function $c: X\times Y \to \Rop, n\in\N$ and $\eps >0$ we set
\begin{equation}\label{eqnbneps}
B_{n,\eps}:=\Bigg\{(x_i,y_i)_{i=1}^n\in (X\times Y)^n : \sum_{i=1}^n c(x_i,y_i) \ge
\sum_{i=1}^n c(x_i,y_{i+1})+\eps\Bigg\}.
\end{equation}
\end{defn}

The definition of the sets $B_{n,\eps}$ is implicitly given in \cite[Theorem 2.3]{GaMc96}. The idea behind it is, that $(x_i,y_i)_{i=1}^n\in B_{n,\eps}$ tells us that transport costs can be reduced if ``$x_i$ is transported to $y_{i+1}$ instead of $y_{i}$'' (recall the conventions $x_{n+1}=x_1$ resp.\ $y_{n+1}=y_1$). In what follows we make this statement precise and give a coordinate free formulation.

Denote by $\sigma, \tau: (X\times Y)^n  \to (X\times Y)^n$ the shifts defined via
\begin{align}
\sigma: (x_i,y_i)_{i=1}^n  &\mapsto  (x_{i+1},y_{i+1})_{i=1}^n \\
  \tau: (x_i,y_i)_{i=1}^n  &\mapsto \hspace{2ex}(x_{i},y_{i+1})_{i=1}^n.
\end{align}
Observe that $\sigma^n=\tau^n=\mbox{Id}_{(X\times Y)^n}$ and that $\sigma$ and $\tau$ commute. Also note that the set $B_{n,\eps}$ from \eqref{eqnbneps} is $\sigma$-invariant (i.e.\ $\sigma(B_{n,\eps})=B_{n,\eps}$), but in general not $\tau$-invariant. Denote by $p_i: (X\times Y)^n \to X\times Y$ the projection on the $i$-th component of the product. The projections $p_X: X\times Y \to X, (x,y)\mapsto x$ and $p_Y: X\times Y \to Y, (x,y)\mapsto y$ are defined as usual and there will be no danger of confusion.

\begin{lem}\label{twoAlternatives}
Let $X,Y$ be Polish  spaces equipped with Borel probability measures $\mu,\nu$. Let $\pi$ be a transport plan. Then one of the following alternatives holds:
\begin{enumerate}
\item $\pi$ is $c$-monotone,

\item there exist $n\in\N$, $\eps>0$ and a measure $\kappa \in \Pi(\pi,\ldots,\pi)$ such that $\kappa(B_{n,\eps})>0$. Moreover $\kappa$ can be taken to be both $\sigma$ and $\tau$ invariant.
\end{enumerate}
\end{lem}
\begin{pf}
Suppose that $B_{n,\eps}$ is an L-shaped null set for all $n\in \N$ and every $\eps>0$. Then there are Borel sets $S_{n,\eps}^1,\ldots, S_{n,\eps}^n\subseteq X\times Y$ of full $\pi$-measure such that \[\left(S_{n,\eps}^1\times\ldots\times S_{n,\eps}^n\right) \cap B_{n,\eps} = \emptyset\] 
 and $\pi $ is concentrated on the $c$-monotone set
\[S=\bigcap_{k=1}^{\infty}\bigcap_{n=1}^\infty\bigcap_{i=1}^n S_{n,1/k}^i.\] 
 If there exist $n\in\N$ and $\eps>0$ such that $B_{n,\eps}$ is not an L-shaped null set, we apply Proposition \ref{PKmeasure} to conclude the existence of a measure $\kappa\in\Pi(\pi,\dotsc,\pi)$ with $\kappa(B_{n,\eps})>0$. To achieve the desired invariance, simply replace $\kappa$  by \begin{equation}\label{eqninv}
\frac{1}{n^2}\sum_{i,j=1}^n (\sigma^i\circ\tau^j)\push\kappa\qed
\end{equation}
\end{pf}
 We are now in the position to prove \textbf{Theorem \ref{CmonEquivOptimal}.a}, i.e.\

\begin{nil} Let $X,Y$ be Polish spaces equipped with Borel probability measures $\mu,\nu$ and let $c: X\times Y \to \Rop$ be a Borel measurable cost function. If $\pi$ is a finite optimal transport plan, then $\pi$ is $c$-monotone.\end{nil}

\begin{pf}
Suppose by contradiction that $\pi$ is optimal, $I_c[\pi]<\infty$ but $\pi$ is not $c$-monotone. Then by Lemma \ref{twoAlternatives} there exist $n\in \N$, $\eps>0$ and an invariant measure $\kappa\in\Pi(\pi,\ldots,\pi)$ which gives mass $\alpha>0$ to the Borel set $B_{n,\eps}\subseteq (X\times Y)^n$. Consider now the restriction of $\kappa$ to $B_{n,\eps}$ defined via $
\hat\kappa(A):= \kappa(A\cap B_{n,\eps})$ for Borel sets $A\subseteq (X\times Y)^n$. $\hat\kappa$ is $\sigma$-invariant since both the measure $\kappa$ and the Borel set $B_{n,\eps}$ are $\sigma$-invariant. Denote the marginal of $\hat\kappa$ in the first coordinate $(X\times Y)$ of $(X\times Y)^n$ by $\hat\pi$. Due to $\sigma$-invariance we have
\[p_i\push\hat\kappa =p_i\push(\sigma\push\hat\kappa)\\
=(p_i\circ \sigma)\push\hat\kappa=p_{i+1}\push\hat\kappa,\]
i.e.\ all marginals coincide and we have $\hat\kappa \in \Pi(\hat\pi,\ldots,\hat\pi)$. Furthermore, since $\hat\kappa\le\kappa$, the same is true for the marginals, i.e.\ $\hat\pi\le \pi$. Denote the marginal of $\tau\push\hat \kappa$ in the first coordinate $(X\times Y)$ of $(X\times Y)^n$ by $\hat\pi_{\beta}$. As $\sigma$ and $\tau$ commute, $\tau\push\hat\kappa$ is  $\sigma$-invariant, so the marginals in the other coordinates coincide with $\hat\pi_{\beta}$. An easy calculation shows that $\hat{\pi}$ and $\hat\pi_{\beta}$ have the same marginals in $X$ resp. $Y$:
\begin{align*}
p_X\push\hat\pi_{\beta}&=p_X\push(p_i\push(\tau\push\hat\kappa))=(p_X\circ p_i\circ\tau)\push \hat\kappa =(p_X\circ p_i)\push \hat\kappa=p_X\push\hat\pi,\\
p_Y\push\hat\pi_{\beta}&=p_Y\push(p_i\push(\tau\push\hat\kappa))=(p_Y\circ p_i\circ\tau)\push \hat\kappa =(p_Y\circ p_{i+1})\push \hat\kappa=p_Y\push\hat\pi.
\end{align*}
The equality of the total masses is proved similarly:
\begin{align*}
\alpha=\hat\pi_{\beta}(X\times Y)&=(p_i\circ\tau)\push \hat\kappa(X\times Y) =p_i\push \hat\kappa(X\times Y)=\hat\pi(X\times Y).
\end{align*}

Next we compute the transport costs associated to $\hat\pi_{\beta}$:

\begin{align*}
\int_{X\times Y} c\, d\hat\pi_{\beta} &= \hspace{5ex}\int_{(X\times Y)^n} c\circ p_1\, d(\tau\push\hat\kappa) & \text{(marginal property)}\\
&= \frac{1}{n}\sum_{i=1}^n \int_{(X\times Y)^n} c\circ p_i\, d(\tau\push\hat\kappa) & \text{($\sigma$-invariance)}\\
&=\frac{1}{n}\sum_{i=1}^n \int_{(X\times Y)^n} (c\circ p_i\circ \tau)\, d\hat\kappa & \text{(push-forward)}\\
&=\frac{1}{n}\sum_{i=1}^n \int_{B_{n,\eps}} (c\circ p_i\circ\tau)\, d\kappa & \text{(definition of $\hat\kappa$)}\\
&\le \frac{1}{n}\int_{B_{n,\eps}}\Big[\sum_{i=1}^n  (c\circ p_i)\,  - \eps\Big]\;d\kappa& \text{(definition of $B_{n,\eps}$)}\\
&= \hspace{5ex}\int_{X\times Y} c d\hat\pi- \eps\frac{\alpha}{n}.&\text{(definition of $\hat\pi$)}
\end{align*}
To improve the transport plan $\pi$ we define
\begin{equation}
\pi_{\beta}:=(\pi-\hat\pi)+\hat\pi_{\beta}.
\end{equation}
Recall that $\pi-\hat\pi$ is a positive measure, so $\pi_{\beta}$ is a positive measure. As $\hat\pi$ and $\hat\pi_{\beta}$ have the same total mass, $\pi_{\beta}$ is a probability measure. Furthermore $\hat\pi$ and $\hat\pi_{\beta}$ have the same marginals, so $\pi_{\beta}$ is indeed a transport plan. We have
\begin{equation}
I_c[\pi_{\beta}]=I_c[\pi]+ \int_{X\times Y} c\,
d(\hat\pi_{\beta}-\hat\pi) \le I_c[\pi]-\eps\frac{\alpha}{n} <
I_c[\pi].\qed
\end{equation}
\end{pf}

\section{Connecting $c$-monotonicity and strong $c$-monotonicity}\label{SmilestoneB}

The Ambrosio-Pratelli example (Example \ref{ExAmPra}) shows that  $c$-monotonicity  need not imply  strong $c$-monotonicity in general. Subsequently we shall present a condition which ensures that this implication is valid.

A $c$-monotone transport plan resists the attempt of enhancement by
means of cyclically rerouting.  This, however, may be due to the fact that cyclical
rerouting is a priori impossible due to infinite transport costs on
certain routes.  

Continuing Villani's interpretation, a situation where rerouting in this consortium of bakeries and caf{\'e}s is possible in a satisfactory way
is as follows: Suppose that bakery $x=x_0$ is able to produce one more croissant than it already does and that caf\'{e} $\tilde y$ is short of one
croissant. It might not be possible to transport the additional
croissant itself to the caf\'e in need, as the costs $c(x, \tilde
y)$ may be infinite. Nevertheless it might be possible to find
another bakery $x_1 $ (which usually supplies caf\'e $y_1$) such that
bakery $x$ can transport (with finite costs!) the extra croissant to
$y_1$; this leaves us with a now unused item from bakery $x_1$,
which can be transported to $\tilde y$ with finite costs. Of course
we allow not only one, but finitely many intermediate pairs
$(x_1,y_1),\dotsc,(x_n,y_n)$ of bakeries/caf\'{e}s to achieve this
relocation of the additional croissant.

In the Ambrosio-Pratelli example we can reroute from a point $(x,x\oplus \alpha)\in \Gamma_1$ to a point $(\tilde x,\tilde x\oplus \alpha)\in \Gamma_1$ only if there exists $n\in \N$ such that $x\oplus (n\alpha)=\tilde x$. In particular, irrationality of $\alpha$ implies that if we can redirect with finite costs from $(x,x\oplus \alpha)$ to $(\tilde x,\tilde x\oplus \alpha)$ we never can redirect back from $(\tilde x,\tilde x\oplus \alpha)$ to $(x,x\oplus \alpha)$. 


\begin{defn}\label{Connectedness}
\standPol,
let $c:X\times Y\to \Rop$ be a Borel measurable cost function
and $\Gamma\subseteq X\times Y$  a Borel measurable set on which
$c$ is finite. We define \begin{enumerate}
\item  $(x,y)\lesssim(\tilde x,\tilde y)$ if there exist pairs
$(x_0,y_0),\dotsc,(x_n,y_n)\in\Gamma$ such that $(x,y)=(x_0,y_0)$
and $(\tilde x,\tilde y)=(x_n,y_n)$ and $c(x_1,y_0),\dotsc,c(x_n,y_{n-1})<\infty$.
\item $(x,y)\approx (\tilde x,\tilde y)$ if $(x,y)\lesssim(\tilde x,\tilde y)$
and $(x,y)\gtrsim (\tilde x,\tilde y)$.
\end{enumerate}
We call $(\Gamma,c)$ \emph{connecting} if $c$ is finite on $\Gamma$ and $(x,
y)\approx(\tilde x,\tilde y)$ for all $(x,y), $ $(\tilde x,\tilde y)\in\Gamma.$
\end{defn}

 These relations were introduced in \cite[Chapter 5, p.75]{Vill05} and appear in a construction due to
Stefano Bianchini.

When there is any danger of confusion we will write $\lesssim_{c,\Gamma}$ and
$\approx_{c,\Gamma}$, indicating the dependence on $\Gamma$ and $c$.  Note that $\lesssim$ is a pre-order, i.e.\ a transitive and reflexive relation, and that $\approx$ is
an equivalence relation. We will also need the projections  $\lesssim_X,\approx_X$ resp.\ $\lesssim_Y, \approx_Y$ of these relations onto the set
$p_X[\Gamma]\subseteq X$ resp.\ $p_Y[\Gamma]\subseteq Y$. The projection is defined in the obvious way: $x\lesssim_X \tilde x$ if there exist $y,\tilde y$ such
that  $(x,y),(\tilde x,\tilde y)\in \Gamma $ and $(x,y)\lesssim (\tilde
x,\tilde y)$ holds.

The other relations are defined analogously. The
projections of $\lesssim$ are again pre-orders and the projections of
$\approx$ are again equivalence relations, provided $c$ is finite on
$\Gamma$. The equivalence classes of $\approx$ and its projections
are compatible in the sense that $[(x,y)]_\approx=
\left([x]_{\approx_X}\!\times [y]_{\approx_Y}\right)\cap \Gamma$. The elementary proofs of these facts are left to the reader.

The main objective of this section is to prove Proposition
\ref{CmonImpliesStrongCmon}, based on several lemmas which will be introduced
throughout the section.

\begin{prop}\label{CmonImpliesStrongCmon} \standPol\ and let $c:X\times Y\to \Rop$ be a Borel measurable cost function. Let $\pi$ be a  finite transport plan. Assume that there exists a $c$-monotone set $\Gamma \subseteq X\times Y$ with $\pi(\Gamma)=1$ on which $c$ is finite, such that  $(\Gamma, c)$ is connecting. Then $\pi$ is strongly $c$-monotone.
\end{prop}

In the proof of Proposition \ref{CmonImpliesStrongCmon} we will establish the existence of the functions $\phi, \psi$ using the construction given in \cite{R96}, see also \cite[Chapter 2]{Vill03} and \cite[Theorem 3.2]{AmPr03}. As we do not impose any continuity assumptions on the cost function $c$, we can not prove the Borel measurability of $\phi$ and $\psi$ by using limiting procedures similar to the methods used in \cite{AmPr03,R96,ScTe08,Vill03}. Instead we will use the following projection theorem, a proof of which can be found in \cite[Theorem III.23]{CaVa77} by analyists or in \cite[Section 29.B]{Kech95} by readers who have some interest in set theory.

\begin{prop}\label{univmeasurable} \hspace{-3mm} \footnote{Sets which are images of Borel sets under measurable functions are called \emph{analytic} in descripitive set theory. Lusin first noticed that analytic sets are universally measurable. Details can be found for instance in \cite{Kech95}.}
Let $X$ and $Y$ be Polish spaces, $A\subseteq X$ a Borel measurable set and $f: X\to Y$ a Borel measurable map. Then $B:=f(A)$ is \emph{universally measurable}, i.e. $B$ is measurable with respect to the completion of every $\sigma$-finite Borel measure on $Y$.\end{prop}


The system of universally measurable sets is a $\sigma$-algebra. If $X$ is a Polish space, we call a function $f:X\to[-\infty,\infty]$ universally measurable if the pre-image of every Borel set is universally measurable.


\begin{lem}\label{borelversion} Let $X$ be a Polish space and $\mu$ a finite Borel measure on $X$.
If $\phi: X\to \Rm$ is universally measurable, then there exists a
Borel measurable function $\tilde{\phi}: X\to\Rm$ such that
$\tilde{\phi} \leq\phi$ everywhere and $\phi=\tilde{\phi}$ almost everywhere.
\end{lem}
\begin{pf}
Let $(I_n)_{n=1}^{\infty}$ be an enumeration of the intervals
$[a,b)$ with endpoints in $\Q$ and denote the completion of $\mu$ by $\tilde \mu$. Then for each $n\in \N$,
$\phi^{-1}[I_n]$ is $\tilde{\mu}$-measurable and hence the union of
a Borel set $B_n$ and a $\tilde{\mu}$-null set $N_n$. Let $N$ be a
Borel null set which covers $\bigcup_{n=1}^{\infty} N_n$. Let
$\tilde{\phi}(x)=\phi(x)-\infty\cdot \mathbf{1}_N(x)$. Clearly
$\tilde{\phi}(x) \le \phi(x)$ for all $x\in X$ and
$\phi(x)=\tilde{\phi}(x)$ for $\tilde\mu$-almost all $x\in X$.
Furthermore, $\tilde{\phi}$ is Borel measurable since
$(I_n)_{n=1}^{\infty}$ is a generator of the Borel $\sigma$-algebra
on $\Rm$ and for each $n\in \N$ we have that
$\tilde{\phi}\me [I_n]=B_n\setminus N$ is a Borel set. \qed\end{pf}

The following definition of the functions $\phi_n, n\in\N$ resp.\ $\phi$ is reminiscent of the construction in \cite{R96}.

\begin{lem}\label{analytic} Let $X,Y$ be Polish spaces, $c: X\times Y \to \Rop$ a Borel measurable cost function and $\Gamma\subseteq X\times Y$ a Borel set. Fix $(x_0,y_0)\in \Gamma$ and assume that $c$ is finite on $\Gamma$. For $n\in\N$, define $\phi_n: X\times\Gamma^{n}\to \Rp $ by
\begin{equation}
\phi_n (x;x_1,y_1,\ldots,x_n,y_n)=
[c(x,y_n)-c(x_n,y_n)]+\sum_{i=0}^{n-1} [
c(x_{i+1},y_i)-c(x_i,y_i)]
\end{equation}
Then the map $\phi: X \to [-\infty, \infty]$ defined by
\begin{equation}\label{phirueschen}
\phi(x)=\inf \{\phi_n (x;x_1,y_1,\ldots,x_n,y_n): n\ge 1, (x_i,y_i)_{i=1}^n \in \Gamma^n\}
\end{equation}
is universally measurable.
\end{lem}
\begin{pf}
First note that the Borel $\sigma$-algebra on $[-\infty, \infty]$
is generated by intervals of the form $[-\infty,\alpha)$, thus it is
sufficient to determine the pre-images of those sets under $\phi$. We
have
\begin{equation*}
\phi(x)< \alpha \leftrightarrow
\exists n\in \N \;\exists (x_1,y_1),\ldots,(x_n,y_n)
\in \Gamma:\phi_n(x;x_1,y_1,\ldots,x_n,y_n)<\alpha.
\end{equation*}
The set $\phi_n^{-1}[[-\infty,\alpha)]$ is Borel measurable. Hence
\[\phi^{-1}[[-\infty,\alpha)]= \bigcup_{n\in\N}
p_X[\phi_n^{-1}[[-\infty,\alpha)]]\] is the countable union of
projections of Borel sets. Since projections of Borel sets are
universally measurable by Proposition \ref{univmeasurable}, $\phi^{-1}[[-\infty,\alpha)]$ belongs also to the
$\sigma$-algebra of universally measurable sets. \qed\end{pf}

\begin{lem}\label{phi finite} Let $X,Y$ be Polish spaces and $c: X\times Y \to \Rop$ a Borel measurable cost function. Suppose $\Gamma$ is $c$-monotone, $c$ is finite on $\Gamma$ and $(\Gamma,c)$ is connecting. Fix $(x_0,y_0)\in\Gamma$. Then the map $\phi$ from \eqref{phirueschen} is finite on $p_X[\Gamma]$. Furthermore
\begin{equation}\label{ctrafoineq}
\phi(x)\le \phi(x')+c(x,y)-c(x',y)\quad \forall x\in X,\,(x',y)\in \Gamma.
\end{equation}
\end{lem}
\begin{pf}
Fix $x\in p_X[\Gamma]$. Since $x_0\lesssim x$ (recall Definition \ref{Connectedness}), we can find
$x_1,y_1,\ldots,x_n,y_n$ such that $\phi_n
(x;x_1,y_1,\ldots,x_n,y_n)<\infty$. Hence $\phi(x)<
\infty$. Proving $\phi(x)>-\infty$ involves some wrestling with notation
but, not very surprisingly, it comes down to applying the fact that
$x\lesssim x_0$. Let $a_1=x$ and choose $b_1,a_2,b_2,\ldots, a_m,b_m$
such that $(a_1,b_1),\ldots,(a_m,b_m)\in \Gamma$ and $c(a_2,b_1),
\ldots,c(a_m,b_{m-1}),c(x,b_m)<\infty$. Assume now that
$x_1,y_1,\ldots,x_n,y_n$ are given such that $\phi_n
(x;x_1,y_1,\ldots,x_n,y_n)<\infty$. Put $x_{n+i}=a_i$ and $y_{n+i}=b_i$
for $i\in\{1,\ldots,m\}$. Due to $c$-monotonicity of $\Gamma$ and the finiteness of all involved terms we have:
\begin{equation*}\label{c monoton}
0\le [c(x_0,y_{n+m})-c(x_{n+m},y_{n+m})]+\sum_{i=0}^{n+m-1}
[c(x_{i+1},y_i)-c(x_i,y_i)],
\end{equation*}
which, after regrouping yields
\begin{multline}\label{eqnhenceandforth}
\alpha:=[c(x_0,b_{m})-c(a_{m},b_{m})]+\sum_{i=1}^{m-1}[c(a_{i+1},b_i)-c(a_i,b_i)]\\
\le [c(x,y_{n})-c(x_{n},y_{n})]+\sum_{i=0}^{n-1}[c(x_{i+1},y_i)-c(x_i,y_i)].
\end{multline}
Note that the right hand side of \eqref{eqnhenceandforth} is just $\phi_n(x;x_1,y_1,\ldots,x_n,y_n)$. Thus passing to the infimum we see that $\phi(x)\ge \alpha>-\infty.$ To prove the remaining inequality, observe that the right hand side of \eqref{ctrafoineq} can be written as
$$\inf\{\phi_n(x;x_1,y_1,\ldots, x_n,y_n):n\ge 1,(x_i,y_i)_{i=1}^n \in \Gamma^n \mbox{ and } (x_n, y_n)=(x',y)\}$$ whereas the left hand side of \eqref{ctrafoineq} is the same, without the restriction $(x_n,y_n)=(x',y)$.
\qed\end{pf}

\begin{lem}\label{analytic2}
Let $X,Y$ be Polish spaces and $c: X\times Y \to \Rop$ a Borel measurable cost function. Let $X_0\subseteq X$ be a non-empty Borel set and let $\phi: X_0\to \R$ be a Borel measurable function. Then the $c$-transform $\psi: Y \to \Rm$, defined as
\begin{equation}
\psi(y):=\inf_{x\in X_0} [c(x,y)-\phi(x)]
\end{equation}
is universally measurable.
\end{lem}

\begin{pf}
As in the proof of Lemma \ref{borelversion} we consider the set
$\psi^{-1}[[-\infty,\alpha)]$:
\[\psi(y)< \alpha \leftrightarrow \exists x\in X_0:
c(x,y)-\phi(x)<\alpha.\] Note that the set $\{(x,y)\in X_0\times Y:
c(x,y)-\phi(x)<\alpha\}$ is Borel. Thus
\[\psi^{-1}[[-\infty,\alpha)]=p_X[\{(x,y)\in X_0\times Y: c(x,y)-\phi(x)<\alpha\}]\] is the projection of a Borel set, hence universally measurable.
\qed\end{pf}

We are now able to prove the main result of this section.

\begin{pf*}{PROOF of Proposition \ref{CmonImpliesStrongCmon}.}
Let $\Gamma\subseteq X\times Y$ be a $c$-monotone Borel set such that $\pi(\Gamma)=1$ and the pair $(\Gamma, c)$ is connecting.
Let $\phi$ be the map from Lemma \ref{analytic}. Using Lemma
\ref{borelversion} and Lemma \ref{phi finite}, and eventually
passing to a subset of full $\pi$-measure, we may assume that
$\phi$ is Borel measurable, that $X_0:=p_X[\Gamma]$ is a Borel set and
that
\begin{equation}\label{equalonC}
c(x',y)-\phi(x')\le c(x,y)-\phi(x)\quad \forall x\in X_0,\,(x',y)\in \Gamma.
\end{equation}
Note that \eqref{equalonC} follows from \eqref{ctrafoineq} in Lemma \ref{phi finite}.  Here we consider $x\in X_0$ in order to ensure that $\phi(x)$ is finite on $X_0$. Now consider the $c$-transform
\begin{equation}\label{eqctrafo}
\psi(y):=\inf_{x\in X_0} [c(x,y)-\phi(x)],
\end{equation}
which by Lemma \ref{analytic2} is universally measurable. Fix $y \in
p_Y[\Gamma]$. Using \eqref{equalonC} we see that the infimum in
\eqref{eqctrafo} is attained at a point $x_0\in X_0$ satisfying
$(x_0,y)\in \Gamma$. This implies that $\phi(x)+\psi(y)=c(x,y)$ on
$\Gamma$ and $\phi(x)+\psi(y)\le c(x,y)$ on $p_X[\Gamma]\times
p_Y[\Gamma]$. To guarantee this inequality on the whole product
$X\times Y$, one has to redefine $\phi$ and $\psi$ to be $-\infty$
on the complement of $p_X[\Gamma]$ resp.\ $p_Y[\Gamma]$. Applying
Lemma \ref{borelversion} once more, we find that there exists a
Borel set $N\subseteq Y$  of zero $\nu$-measure, such that
$\tilde\psi(y)=\psi(y)-\infty\cdot \Eins_N(y)$ is Borel measurable.
Finally, replace $\Gamma$ by $\Gamma \cap (X\times (Y\setminus N))$
and $\psi$ by $\tilde\psi$. \qed\end{pf*}

We conclude this section by proving that every strongly $c$-monotone transport plan is optimal ({\bf Proposition \ref{StrongCmonImpliesOptimal}}). 

\begin{nil}
\standPol\ and
let $c: X\times Y\to \Rop$ be Borel measurable. Then every finite
transport plan which is strongly $c$-monotone is optimal.
\end{nil}

\begin{pf}
Let $\pi_0$ be a strongly $c$-monotone transport plan.  Then, according to the definition, there exist Borel functions $\phi(x)$ and $\psi(y)$ taking values in $[-\infty,\infty)$
such that
\begin{equation}\label{defstrongc}
\phi(x)+\psi(y)\le c(x,y)
\end{equation}
everywhere on $X\times Y$  and equality holds $\pi_0$-a.e.
We define the truncations $\phi_n = (n\wedge (\phi\vee -n)), \psi_n=(n\wedge (\psi\vee -n))$ and let $\xi_n(x,y):=\phi_n(x)+\psi_n(y)$ resp.\ $\xi(x,y):=\phi(x)+\psi(y)$. Note that $\phi_n,\psi_n,\xi_n, \xi$ are Borel measurable.
By elementary considerations which are left the reader, we get pointwise monotone convergence $\xi_n\uparrow \xi$ on the set $\{\xi\ge 0\}$ resp.\ $\xi_n\downarrow \xi$ on the set $\{\xi\le 0\}$
Let $\pi_1$ be an arbitrary finite transport plan; to compare $I_c[\pi_0]$ and $I_c[\pi_1]$ we make the following observations:
\begin{enumerate}
\item By monotone convergence
\begin{align}
\int_{\{\xi\geq 0\}} \xi_n \, d\pi_i\uparrow & \int_{\{\xi\geq 0\}} \xi \, d\pi_i \leq I_c[\pi_i]< \infty \mbox{ and}\\
\int_{\{\xi< 0\}} \xi_n \, d\pi_i\downarrow & \int_{\{\xi< 0\}} \xi \, d\pi_i  
\end{align}
for $i\in \{0,1\}$, hence  $\lim_{n\to \infty} \int \xi_n\, d\pi_i= \int \xi\, d\pi_i.$

\item By the assumption on equal marginals of $\pi_0$
and $\pi_1$ we obtain for $n\geq 0$
\begin{align}
\int \xi_n\,d\pi_0  &  =\int \phi_{n}\, d\pi_0+\int \psi_{n}\,d\pi_0\\
&  =\int \phi_{n}\, d\pi_1+\int \psi_{n}\,d\pi_1  =\int\xi_{n}\,d\pi_1.
\end{align}
\end{enumerate}
 Thus $I_c[\pi_0]=\int \xi\, d\pi_0= \lim_{n\to \infty} \int \xi_n\, d\pi_0= \lim_{n\to \infty} \int \xi_n\, d\pi_1=\int \xi\, d\pi_1\leq I_c[\pi_1]$; since $\pi_1$ was arbitrary,  this implies optimality of $\pi_0$.
\qed\end{pf}

\section{From $c$-monotonicity to optimality}\label{SmilestoneC}

This section is devoted to the proof of Theorem \ref{CmonEquivOptimal}.b.
Our argument starts with a finite $c$-monotone transport plan $\pi$
and we aim for showing that $\pi$ is at least as good as any other
finite transport plan. The idea behind the proof is to partition $X$
and $Y$ into cells $C_i, i\in I$ resp.\ $D_i,i\in I$ in such a way
that $\pi$ is strongly  $c$-monotone on ``diagonal'' sets of the
form $C_i\times D_i$ while regions $C_i\times D_j, i\neq j$ can be
ignored, because no finite transport plan will give positive measure
to the set $C_i\times D_j$.

Thus it will  be necessary to apply previously established results
to some restricted transport problems on a space $C_i\times D_i$
equipped with some relativized transport plan $\pi\restr C_i\times
D_i$. As in general  the cells $C_i,D_i$ are  plainly Borel sets they
may fail to be Polish spaces with respect to the topologies
inherited from $X$ resp.\ $Y$. However, for us it is only important
that there exist \emph{some} Polish topologies that generate the
same Borel sets on $C_i$ resp.\ $D_i$ (see e.g.\ \cite[Theorem 13.1]{Kech95}). At this point it is crucial that our results only
need measurability of the cost function and do not ask for any form
of continuity (cf.\ the remarks at the end of the introduction).
Before we give the proof of Theorem \ref{CmonEquivOptimal}.b we will
need some preliminary lemmas.

\begin{lem}\label{omegaClasses} \standPol\ and let $c: X\times Y \to \Rop$ be a Borel measurable cost function.
Let $\pi,\pi_0$ be finite transport plans and $\Gamma\subseteq X\times Y$ a Borel set with $\pi(\Gamma)=1$ on which $c$ is finite.
Let $I=\{0,\ldots, n\}$ or $I=\N$
and assume that $C_i, i\in I $ are mutually disjoint Borel sets in $X$, $D_i,i\in  I$ are mutually disjoint Borel sets in $Y$ such that the equivalence classes of $\approx_{c,\Gamma}$ are of the form $\Gamma \cap (C_i\times D_i)$.
Then also $\pi_0(\Cup_{i\in I } C_i\times D_i)=1$.
\end{lem}

In the proof we will need the following simple lemma. (For a proof see for instance \cite[Proposition 8.13]{Kall97}.)
\begin{lem}\label{markoversatz}
Let $I=\{0,\ldots,n\}$ or $I=\N$ and let $P=(p_{ij})_{i,j\in I} $ be a matrix with non-negative entries such that $\sum_{j\in I} p_{i_0j}=1$ for each $i_0\in I$. 
Assume that there exists a vector $(p_i)_{i\in I}$ with strictly positive entries such $p\cdot P= p$.\footnote{Such a matrix $P$ is often called  a \emph{stochastic matrix} while $p$ is a \emph{stochastic vector}.}  
Then whenever $p_{i_0 i_1}>0$ for $i_0, i_1\in I$, there
exists a finite sequence $i_0, i_1, \ldots, i_n=i_0$ such that for all $0\le k < n$ one has $p_{i_k i_{k+1}}>0$.
\end{lem}

\begin{pf*}{PROOF of Lemma \ref{omegaClasses}.}
As $\approx_{\Gamma, c}$ is an equivalence relation and $\pi$ is concentrated on $\Gamma$, the sets $C_i$, $i\in I$ are a partition of $X$ modulo $\mu$-null sets. Likewise the sets  $D_i$, $i\in I$ form a partition of $Y$ modulo $\nu$-null sets.
In particular the quantities
\begin{equation}
p_i:=\mu(C_i)=\nu(D_i)=\pi (C_i\times D_i),\quad i\in I
\end{equation}
add up to $1$. Without loss of generality we may assume that
$p_i>0$ for all $i\in I$. We define
\begin{equation}
p_{ij}:=\frac{\pi_0(C_i\times D_j)}{\mu(C_i)}, \quad  i,j\in I.
\end{equation}
Then $\sum_{j\in I} p_{i_0j}= \frac{\pi_0(C_i\times Y)}{\mu(C_i)}=1$ for each $i_0\in I$.  By the condition on the marginals of
$\pi_0$ we have for the $i$-th component of $p\!\cdot\!P$
\[(p\!\cdot\!P)_i = \sum_{j\in I}\mu(C_j) \frac{\pi_0(C_j\times D_i)}{\mu(C_j)}= \pi_0(X\times D_i) = \nu (D_i)=p_i\]
i.e.\ $p\!\cdot\!P=p$. Hence $P$ satisfies the assumptions of Lemma \ref{markoversatz}. 
We claim that $p_{ii}=1$ for all $i\in I$.
 Suppose not. Pick $i_0\in I$ such that $p_{i_0i_0}<1$.  Then there exists some index $i_1\neq i_0$ such that $p_{i_0i_1}>0$. 
 Pick a finite sequence $i_0,i_1,\ldots, i_n=i_0$ according to Lemma \ref{markoversatz}. Fix $k\in \{1,\dotsc, n-1\}$. Then $$\pi_0(C_{i_k}\times D_{i_{k+1}})=p_{i_k i_{k+1}}>0.$$
Since $\pi_0$ is a finite transport plan, there exist $x_k\in C_{i_k}\cap p_X[\Gamma]$ and $y_{k+1}'\in D_{i_{k+1}}\cap p_Y[\Gamma]$ such that $c(x_k,y_{k+1}')<\infty$. Choose $y_k\in D_{i_k}$ and $x_{k+1}'\in C_{i_{k+1}}$ such that $(x_k,y_k), (x_{k+1}',y_{k+1}')\in \Gamma.$ Then
 \begin{equation*}(x_0,y_0)\lesssim (x_1',y_1')\approx(x_1,y_1)\lesssim (x_2',y_2')\approx(x_2,y_2)\lesssim\dotsc\lesssim (x_n',y_n')\approx (x_0,y_0).
 \end{equation*}
 But this implies that $(x_0,y_0)\approx (x_1,y_1),$ contradicting the assumption that  $(C_{i_0}\times D_{i_0})\cap \Gamma, (C_{i_1}\times D_{i_1})\cap \Gamma$ are different equivalence classes of $\approx_{\Gamma,c}$. Hence we have indeed $p_{ii}=1$ for all $i\in I$, thus $ \pi_0(C_i\times D_i)=\mu(C_i)$ which implies $\pi_0(\bigcup_{i\in I} C_i\times D_i)=1$.
\qed\end{pf*}
\begin{lem}\label{CFiniteImpliesConnected}
\standPol\ and let $c:X\times Y\to [0,\infty]$ be a Borel measurable cost function which is $\mu\otimes\nu$-a.e.\ finite. For every finite transport plan $\pi$ and every Borel set $\Gamma\subseteq X\times Y$ with $\pi(\Gamma)=1$ on which $c$ is finite, there exist Borel sets $O\subseteq X,U\subseteq Y$ such that  $\Gamma'=\Gamma\cap (O\times U)$ has full $\pi$-measure and $(\Gamma',c)$ is connecting. 
\end{lem}

\begin{pf} 
By Fubini's Theorem for $\mu$-almost all $x\in X$ the set $\{y: c(x,y)<\infty\}$ has full $\nu$-measure and for $\nu$-almost all $y\in Y$ the set $\{x: c(x,y)<\infty\}$ has full $\mu$-measure. In particular the set of points $(x_0,y_0)$ such that both $\mu\left(\{x:
c(x,y_0)<\infty\}\right)=1$ and $\nu\left(\{y: c(x_0,y)<\infty\}\right)=1$ has full
$\pi$-measure. Fix such a pair $(x_0,y_0)\in\Gamma$ and let
$O=\{x\in X: c(x,y_0)<\infty  \}, U=\{y\in Y:c(x_0,y)<\infty\}$.
 Then $\Gamma'=\Gamma \cap (O\times U)$ has full $\pi$-measure and for
every $(x,y)\in\Gamma'$ both quantities $c(x,y_0)$ and $c(x_0,y)$
are finite. Hence $x\approx_X x_0$, for every $x\in p_X[\Gamma']$.
Similarly we obtain $y\approx_Y y_0$, for every $y\in p_Y[\Gamma']$.
Hence $(\Gamma',c)$ is connecting. \qed\end{pf}

Finally we  prove the statement of \textbf{Theorem \ref{CmonEquivOptimal}.b}:

\begin{nil}
\standPol\ and $c: X\times Y \to \Rop$ a Borel measurable cost function. Every finite $c$-monotone transport plan is optimal if there exist a closed set $F$ and a $\mu\otimes\nu$-null set $N$ such that $\{(x,y): c(x,y)=\infty\} = F\cup N$.
\end{nil}

\begin{pf}
Let $\pi$ be a finite $c$-monotone transport plan and pick a $c$-monotone Borel set $\Gamma\subseteq X\times Y$ with $\pi(\Gamma)=1$
on which $c$ is finite.

Let $O_n,U_n,n\in \N$ be open sets such that $\bigcup_{n\in \N} (O_n\times U_n)=(X\times Y)\setminus F$.
Fix $n\in \N$ and
interpret $\pi\restr O_n\times U_n$ as a transport plan on the
spaces $(O_n,\mu_n)$ and $(U_n,\nu_n)$ where $\mu_n$ and $\nu_n$ are
the marginals corresponding to $\pi\restr O_n\times U_n$. Apply
Lemma \ref{CFiniteImpliesConnected} to $\Gamma\cap(O_n\times U_n)$ and  the cost
function $c\restr O_n\times U_n$ to find $O_n'\subseteq
O_n,U_n'\subseteq U_n$ and $\Gamma_n =\Gamma\cap (O_n' \times U_n')$
with $\pi(\Gamma_n)=\pi(\Gamma\cap (O_n \times U_n))$ such that
$(\Gamma_n,c)$ is connecting. Then $\tilde
\Gamma=\Cup_{n\in\N}\Gamma_n $ is a  subset of $\Gamma$ of full
measure and every equivalence class of $\approx_{\tilde \Gamma,c} $
can be written in the form $((\bigcup_{n\in N} O_n')\times
(\bigcup_{n\in N} U_n'))\cap \Gamma$ for some non-empty index set
$N\subseteq \N$. Thus there are at most countably many equivalence
classes which we can write in the form $(C_i\times D_i) \cap \Gamma,
i\in I $ where $I=\{1,\dotsc, n\}$ or $I=\N$. Note that by shrinking
the sets $C_i, D_i, i\in I$  we can assume that $C_i\cap C_j=D_i\cap
D_j=\emptyset$ for $i\neq j$.

Assume now that we are given another finite transport plan $ \pi_0$.
Apply Lemma \ref{omegaClasses} to $\pi, \pi_0$ and $\tilde \Gamma$
to achieve that $\pi_0$ is concentrated on $\Cup_{i\in I}
C_i\times D_i$. For $i\in I$ we consider the restricted problem of
transporting $\mu\!\restr\! C_i$ to $\nu\!\restr\!D_i$. We know that
$\pi\! \restr\! C_i\times D_i $ is optimal for this task by
Propositions \ref{StrongCmonImpliesOptimal} and \ref{CmonImpliesStrongCmon}, hence $I_c[\pi]\le I_c[\pi_0]$.
\qed\end{pf}

\begin{rem} In fact the following somewhat more general (but also more complicated to state) result holds true: \emph{Assume that $\{(x,y): c(x,y)=\infty\} \subseteq  F\cup N$ where $F$ is closed and $N$ is a $\mu\otimes\nu$-null set. Then every $c$-monotone transport plan $\pi$ with $\pi(F\cup N)=0$ is optimal.}
\end{rem}

\section{Completing the picture}\label{Spuzzle}
First we give the proof of  \textbf{Theorem
\ref{StrongCmonEquivRobustlyOptimal}}.

\begin{nil}
\standPol\ and $c:X\times Y\to \Rop$ a
Borel measurable cost function. For a finite transport plan $\pi$
the following assertions are equivalent:
\begin{enumerate} 
\item $\pi$ is robustly optimal.
\item $\pi$ is strongly $c$-monotone.
\end{enumerate}
\end{nil}

\begin{pf}
$a. \Rightarrow b.$: Let $Z$ and $\lambda\neq 0$ be according to the definition of robust optimality. As $\tilde\pi=(\mbox{Id}_Z\times\mbox{Id}_Z)\push \lambda + \pi$ is optimal, Theorem \ref{CmonEquivOptimal}.a ensures the existence of a $\tilde c$-monotone Borel set $\tilde \Gamma\subseteq (X\cup Z) \times (Y\cup Z)$ such that $\tilde c$ is finite on $\tilde \Gamma $ and $\tilde \pi $ is concentrated on $\tilde \Gamma$. Note that $(z,z)\in \tilde\Gamma$ for $\lambda$-a.e.\ $z\in Z$.   We claim that for $\lambda$-a.e.\ $z\in Z$ and all $(x,y)\in \Gamma=\tilde \Gamma\cap (X\times Y)$ the relation
\begin{align}
(x,y)\approx_{\tilde \Gamma, \tilde c} (z,z)\end{align} 
holds true. 
Indeed, since $\tilde c$ is finite on $ Z\times Y$ we have $c(z,y)<\infty$ hence $(x,y)\lesssim_{\tilde \Gamma, \tilde c} (z,z)$. Analogously finiteness of $\tilde c$ on $X\times Z$ implies $c(x,z)<\infty$ such that also $(z,z)\lesssim_{\tilde \Gamma, \tilde c} (x,y)$.
By transitivity of $\approx_{\tilde \Gamma, \tilde c}$, $(\tilde
\Gamma,\tilde c)$ is connecting. Applying Proposition
\ref{CmonImpliesStrongCmon} to the spaces $X\cup Z$ and $Y\cup Z$ we get that
$\tilde\pi$ is strongly $\tilde c$-monotone, i.e.\ there exist
$\tilde \phi$ and $\tilde \psi$ such that $\tilde\phi(a)+\tilde\psi(b)\le \tilde c(a,b)$ for $(a,b)\in (X\cup Z)\times (Y\cup Z)$  and equality holds
$\tilde \pi$-almost everywhere. By restricting $\tilde \phi$ and
$\tilde \psi$ to $X$ resp.\ $Y$ we see that $\pi$ is strongly
$c$-monotone.
\medskip

\noindent
$b. \Rightarrow a.$: Let $Z$ be a Polish space and let $\lambda$ be a  finite Borel measure on $Z$.  We extend $c$ to $\tilde c: (X\cup Z)\times (Y\cup Z)\to \Rop$ via
\[\tilde{c}(a,b)=\left\{\begin{array}{cl}
c(a,b)                    &\mbox{ for } (a,b)\in X\times Y\\
\max\left(\phi(a),0\right)&\mbox{ for } (a,b)\in X\times Z\\
\max\left(\psi(b),0\right)&\mbox{ for } (a,b)\in Z\times Y\\
0     &\mbox{ otherwise.}  \end{array} \right.\]

Define $\tilde\phi(a):=\left\{\begin{array}{cll}
\phi(a)&\mbox{ for }& a\in X\\
0      &\mbox{ for }& a \in Z\end{array} \right.$ and
$\tilde\psi(b):=\left\{\begin{array}{cll}
\psi(b)&\mbox{ for }& b\in Y\\
0     &\mbox{ for }&b\in Z.\end{array} \right.$

Then $\tilde\phi$ resp.\ $\tilde\psi$ are extensions of $\phi$ resp.
$\psi$ to $X\cup Z$ resp. $Y\cup Z$ which satisfy
$\tilde\phi(a)+\tilde\psi(b)\le \tilde c(a,b)$ and equality holds on
$\tilde\Gamma=\Gamma\cup\{(z,z):z\in Z\}$. Hence $\tilde \Gamma$ is
strongly $\tilde c$-monotone. Since $\tilde{\pi}$ is concentrated on
$\tilde{\Gamma}$, $\tilde \pi $ is optimal by Proposition \ref{StrongCmonImpliesOptimal}. \qed\end{pf}

Next consider \textbf{Theorem \ref{AllEquiv}}.

\begin{nil}
\standPol\ and let $c: X\times Y\to \Rop$ be Borel measurable and
$\mu\otimes\nu$-a.e.\ finite. For a finite transport plan $\pi$ the
following assertions are equivalent:
\begin{enumerate}\renewcommand{\labelenumi}{(\arabic{enumi})}
\item $\pi$ is optimal.
\item $\pi$ is $c$-monotone.
\item $\pi$ is robustly optimal.
\item $\pi$ is strongly $c$-monotone.
\end{enumerate}
\end{nil}
\begin{pf}  
By Theorem \ref{StrongCmonEquivRobustlyOptimal}, (3) and (4) are equivalent and they trivially  imply (1) and (2) which are equivalent by Theorem \ref{CmonEquivOptimal}. It remains to see that $(2) \Rightarrow (4)$. Let $\pi $ be a finite $c$-monotone transport plan. Pick a $c$-monotone Borel set $\Gamma \subseteq X\times Y$ such that $c$ is finite on $\Gamma $  and $\pi(\Gamma)=1$. By Lemma \ref{CFiniteImpliesConnected} there exists a Borel set $\Gamma'\subseteq \Gamma$ such that $\pi(\Gamma')=1$ and $(\Gamma',c)$ is connecting, hence Proposition \ref{CmonImpliesStrongCmon} applies.
\qed\end{pf}

Finally the example below shows that the ($\mu\!\otimes\!\nu$-a.e.) finiteness of the cost function is essential to be able to pass from the ``weak properties'' (optimality, $c$-monotonicity) to the ``strong properties'' (robust optimality, strong $c$-monotonicity).

\begin{exmp}[Optimality does not imply strong $c$-monotonicity]
\label{ZeroOneExample} Let $X=Y=[0,1]$ and equip both spaces with
Lebesgue measure $\lambda=\mu=\nu$. Define $c$ to be $\infty $ above
the diagonal and $1-\sqrt{x-y}$ for $y\leq x$. The optimal (in this
case the only finite) transport plan is the Lebesgue measure $\pi$
on the diagonal $\Delta$. We claim that $\pi$ is not strongly
$c$-monotone. Striving for a contradiction we assume that there
exist $\phi$ and $\psi$ witnessing the strong $c$-monotonicity. Let
$\Delta_1$ be the full-measure subset of $\Delta$ on which
$\phi+\psi=c$, and write $p_X[\Delta_1]$ for the
 projection of $\Delta_1$.
We claim that
\begin{equation}\label{localineq}
\forall x,x'\in p_X[\Delta_1]: \mbox{ If } x<x', \mbox{ then }
\phi(x)-\phi(x') \geq \sqrt {x'-x},
\end{equation}
which will yield a contradiction when combined with the fact that $p_X[\Delta_1]$
is dense.

Our claim  \eqref{localineq} follows directly  from
\begin{equation}
\phi(x')+\psi(x)\leq c(x',x)= 1-\sqrt{x'-x}
 \ \mbox{ and }\ \phi(x)+\psi(x) = c(x,x)=1.
 \end{equation}

Now let $x<x+a $ be elements of $p_X[\Delta_1]$, let $b:= \phi(x)-\phi(x')$, and
let $n\in \N$ be a sufficiently large number, say
satisfying $n>2\tfrac{b^2}{a^2}$.  Using the fact
that $p_X[\Delta_1]$ is dense, we can find real numbers
$x=x_0< x_1 < \cdots <  x_n = x+a$
in $\Delta_1$ satisfying $x_k-x_{k-1} < 2/n$ for $k=1,\ldots, n$.

Let $\eps_k:= x_k-x_{k-1}$ for $k=1,\ldots, n$. Then we have $\eps_k < \tfrac 2n <
\tfrac {a^2}{b^2}$ for all $k$, hence  $ \sqrt{\eps_k} > \tfrac ba \eps_k$.
So we get
\begin{equation*}
b = \phi(x)-\phi(x') = \sum_{k=1}^n \phi(x_{k-1})-\phi(x_{k})  \ge
\sum_{k=1}^n \sqrt{\eps_k}   >
\sum_{k=1}^n  \frac{b}{a}  \eps_k  =
\frac{b}{a} \sum_{k=1}^n \eps_k = b,
\end{equation*}
a contradiction. (By letting $c=0$ below the  diagonal the argument could be simplified, but
then we would lose lower semi-continuity of $c$.)
\end{exmp}

\def\ocirc#1{\ifmmode\setbox0=\hbox{$#1$}\dimen0=\ht0 \advance\dimen0
  by1pt\rlap{\hbox to\wd0{\hss\raise\dimen0
  \hbox{\hskip.2em$\scriptscriptstyle\circ$}\hss}}#1\else {\accent"17 #1}\fi}


\begin{thebibliography}{10}
\expandafter\ifx\csname url\endcsname\relax
  \def\url#1{\texttt{#1}}\fi
\expandafter\ifx\csname urlprefix\endcsname\relax\def\urlprefix{URL }\fi

\bibitem{AmPr03}
L.~Ambrosio, A.~Pratelli, Existence and stability results in the {$L\sp 1$}
  theory of optimal transportation, in: Optimal transportation and applications
  (Martina Franca, 2001), vol. 1813 of Lecture Notes in Math., Springer,
  Berlin, 2003, pp. 123--160.

\bibitem{CaVa77}
C.~Castaing, M.~Valadier, Convex analysis and measurable multifunctions,
  Lecture Notes in Mathematics, Vol. 580, Springer-Verlag, Berlin, 1977.

\bibitem{GaMc96}
W.~Gangbo, R.~J. McCann, The geometry of optimal transportation, Acta Math.
  177~(2) (1996) 113--161.

\bibitem{Kall97}
O.~Kallenberg, Foundations of modern probability, Probability and its
  Applications (New York), Springer-Verlag, New York, 1997.

\bibitem{Kech95}
A.~S. Kechris, Classical descriptive set theory, vol. 156 of Graduate Texts in
  Mathematics, Springer-Verlag, New York, 1995.

\bibitem{Kell84}
H.~G. Kellerer, Duality theorems for marginal problems, Z. Wahrsch. Verw.
  Gebiete 67~(4) (1984) 399--432.

\bibitem{KnSm92}
M.~Knott, C.~Smith, On {H}oeffding-{F}r\'echet bounds and cyclic monotone
  relations, J. Multivariate Anal. 40~(2) (1992) 328--334.

\bibitem{Prat07}
A.~Pratelli, On the sufficiency of $c$-cyclical monotonicity for optimality of
  transport plans, Math. Z. In press.

\bibitem{RaRu98}
S.~T. Rachev, L.~R{\"u}schendorf, Mass transportation problems. {V}ol. {I},
  Probability and its Applications (New York), Springer-Verlag, New York, 1998.

\bibitem{Rama02}
D.~Ramachandran, Perfect measures and related topics, in: Handbook of measure
  theory, Vol. I, II, North-Holland, Amsterdam, 2002, pp. 765--786.

\bibitem{Rock66}
R.~T. Rockafellar, Characterization of the subdifferentials of convex
  functions, Pacific J. Math. 17 (1966) 497--510.

\bibitem{R96}
L.~R{\"u}schendorf, On {$c$}-optimal random variables, Statist. Probab. Lett.
  27~(3) (1996) 267--270.

\bibitem{ScTe08}
W.~Schachermayer, J.~Teichmann, Characterization of optimal transport plans for
  the {M}onge-{K}antorovich problem, Proceedings of the A.M.S. In press.

\bibitem{Vill03}
C.~Villani, Topics in optimal transportation, vol.~58 of Graduate Studies in
  Mathematics, American Mathematical Society, Providence, RI, 2003.

\bibitem{Vill05}
C.~Villani, Optimal Transport Old and New, Series: Grundlehren der mathematischen Wissenschaften, Vol 338, Springer Verlag, Berlin, 2009.

\end{thebibliography}
\end{document}